\documentclass{article}
\usepackage{amsmath, epsfig,amssymb, multicol, amsthm}
\newtheorem{lemma}{Lemma}
\newtheorem{theorem}{Theorem}

\newcommand{\E}{{\rm E}}

\renewcommand{\Pr}{{\rm Pr}}
\newcommand{\Po}{{\rm Po}}
\newcommand{\D}{{\cal D}}

\renewcommand{\l}{{\cal L}}

\newcommand{\N}{{\cal N}}
\newcommand{\R}{{\cal R}}

\renewcommand{\P}{{\cal P}}

\newcommand{\X}{{\cal X}}

\title{High-order phase transition in random hypergraphs}
\author{Linyuan Lu
\thanks{Department of Mathematics, University of South Carolina, Columbia, SC 29208, USA
({\tt lu@math.sc.edu}). This author was supported in part by NSF
grant DMS 1600811.}
  \and Xing Peng
\thanks{Center for Applied Mathematics, Tianjin University,  Tianjin 300072, P.~R.~China,
({\tt x2peng@tju.edu.cn})Research is supported in part by National Natural Science Foundation of China (No. 11601380). }
}
\date{}
\begin{document}
\maketitle

\begin{abstract}
In this paper, we study the high-order phase transition in random
$r$-uniform hypergraphs. For a positive integer $n$ and a real $p\in [0,1]$,
let $H:=H^r(n,p)$ be the random $r$-uniform  hypergraph
with vertex set $[n]$, where each $r$-set is selected as an edge with
probability $p$ independently randomly. For $1\leq s \leq r-1$
and two $s$-sets $S$ and $S'$, we say $S$  is connected to  $S'$  if
there is a sequence of alternating $s$-sets and edges $S_0,F_1,S_1,F_2, \ldots, F_k,
S_k$ such that $S_0,S_1,\ldots, S_k$ are $s$-sets, $S_0=S$, $S_k=S'$,
 $F_1,F_2,\ldots, F_k$ are edges of $H$, and $S_{i-1}\cup S_i\subseteq
 F_i$ for each  $1\leq i\leq k$. This is an equivalence relation over
the family of all $s$-sets ${[n]\choose s}$ and results in  a partition:
${V\choose s}=\cup_i C_i$. Each $C_i$ is called an { $s$-th-order} connected
component and a component $C_i$ is {\em giant}  if $|C_i|=\Theta(n^s)$. We prove that the sharp threshold of
the existence of the $s$-th-order giant connected
components in $H^r(n,p)$ is $\frac{1}{\big({r\choose s}-1\big){n\choose
    r-s}}$. Let $c={n\choose r-s}p$.
If $c$ is a constant and $c<\tfrac{1}{\binom{r}{s}-1}$, then with high probability,
all $s$-th-order connected components have size $O(\ln n)$.
If $c$ is a constant and $c > \tfrac{1}{\binom{r}{s}-1}$, then with high probability,  $H^r(n,p)$ has a unique giant connected
$s$-th-order component and its size is $(z+o(1)){n\choose s}$, where
$$z=1-\sum_{j=0}^\infty \frac{\left({r\choose
    s}j -j+1 \right)^{j-1}}{j!}c^je^{-c\left({r\choose s}j
  -j+1\right)}.$$
\end{abstract}

\section{Introduction}
The theory of random graphs was born when Erd\H{o}s and R\'enyi
wrote a series of remarkable papers \cite{ER1,ER2}   on the evolution of random graphs
around 1960. In their original papers,  Erd\H{o}s and R\'enyi considered
the uniform model $G_{n,m}$ where a random graph $G$ is selected
uniformly among all graphs with $n$ vertices and $m$
edges. Later, the binomial model $G(n,p)$ became the de facto random
graph model, which is also referred as Erd\H{o}s--R\'enyi random graph model.
In $G(n,p)$, each pair of vertices becomes an edge with
probability $p$ independently. It is remarkable that the
random graph $G(n,p)$ experiences the phase transition as $p$ passes through the threshold
$\frac{1}{n}$. When $p<\tfrac{1-\epsilon}{n}$, with high probability all connected components
of $G(n,p)$ are of order $O(\ln n)$; when $p = \tfrac{1}{n}$, with high probability the largest connected component has size $\Theta(n^{2/3})$ \cite{bo,JKLP,L};
when $p>\tfrac{1+\epsilon}{n}$, with high probability
 there is a unique giant component of size $\Theta(n)$ while all
other connected components have size $O(\ln n)$. For a more detailed description of the phase transition phenomenon in $G(n,p)$, see Chapter 11 of  \cite{AS}.

There are a lot of literature on the phase transition of other random
graphs. Aiello, Chung, and Lu \cite{AFL} introduced a general random
graph model--the random graph with given expected degree sequence.   Chung and
Lu \cite{CL1} studied the connected components in this random graph model and determined the size of the giant
connected component in the supercritical phase
\cite{CL2}. Bollob\'as, Janson, and  Riordan \cite{BJR} investigated
the phase transition phenomenon in inhomogeneous random graphs with a
given kernel function.

There are some attempts at generalizing  Erd\H{o}s and R\'enyi's
seminal work on random graphs to random hypergraphs. Let $H^r(n,m)$ be the
random $r$-uniform hypergraph such that each $r$-uniform hypergraph
with $n$ vertices and $m$ edges is selected randomly with the same probability.
The phase transition of the random hypergraph $H^r(n,m)$ was first
analyzed by Schmidt-Pruzan and Shamir \cite{PS}. Namely, they proved
that for $m< \tfrac{n}{r(r-1)}$, with high probability the largest connected component
in $H^r(n,m)$ is of size $O(\ln n)$; for $m = \tfrac{n}{r(r-1)}$, with high probability the largest connected component is of order $\Theta(n^{2/3})$; for
$m>\tfrac{n}{r(r-1)}$, with high probability there is a unique giant connected
component with $\Theta(n)$ vertices. Karo\'nski and {\L}uczak
\cite{KL} took a closer look at the connected components of $H^r(n,m)$
when $m$ is near $\tfrac{n}{r(r-1)}$. The connected components in
 results stated above  are in the following sense.  Let
$H$ be a uniform hypergraph, two vertices $u$ and $u'$ are {\em
connected} if there is an (alternating) sequence of vertices and edges
$v_0, F_1, v_1, F_2, \ldots, F_k, v_k$ such that $u=v_0$, $u'=v_k$, and
$\{v_{i-1},v_i\}\subseteq F_i\in E(H)$ for $1 \leq i \leq k$.

In a hypergraph $H=(V,E)$, high-order connections exist besides the
vertex-to-vertex connection.
For $1\leq s \leq r-1$ and two $s$-sets
$S$ and $S'$, we say $S$ and $S'$ are {\em connected} if there is a sequence of
alternating $s$-sets and edges $S_0,F_1,S_1,F_2, \ldots, F_k, S_k$
such that $S_0,S_1,\ldots, S_k$ are $s$-sets, $S_0=S$, $S_k=S'$;
$F_1,F_2,\ldots, F_k$ are edges of $H$, and $S_{i-1}\cup S_i\subseteq
F_i$ for $1\leq i\leq k$. It is easy to verify this is an equivalence
relation over the family of all $s$-sets ${V\choose s}$ and then
it results in a partition ${V\choose s}=\cup_i C_i$. Each $C_i$ is called an
{\em $s$-th-order connected component} and each $s$-set is called a {\it stop}. We say a component $C_i$ is a
{\em giant $s$-th-order} component if $|C_i|=\Theta(n^s)$.
Alternately, we construct a auxiliary  graph $H^s$ with the
vertex set ${V\choose s}$ and  a pair of two $s$-sets
$\{S, S'\}$ forms an edge in $H^s$
if $S\cup S'\subseteq F$ for some $ F \in E(H)$.
The $s$-th-order connected
components of $H$ are exactly the connected components of the auxiliary graph
$H^s$. A path (or tree) in $H^s$ is called an {\em $s$-path} (or {\em $s$-tree}) of $H$
correspondingly.

In this paper, we consider  the ``binomial model'' of random $r$-uniform
hypergraphs $H^r(n,p)$.
 Namely, $H^r(n,p)$ has the vertex set  $[n]$ and each $r$-set of $[n]$ becomes an edge
   with probability $p$ independently randomly.
 For additional information on random hypergraphs, the reader is referred to the survey  wrote by Karo{\'n}ski and  {\L}{uczak}\cite{KL0}.
 Our goal is to study the
phase transition of the $s$-th-order connected components in $H^r(n,p)$
as $p$ increases from $0$ to $1$.   Bollob\'as and Riordan  (see page 442 of \cite{BR2}) claimed the following.  The branching process arguments in \cite{BR1} can show that the threshold for thee emergence of the $s$-th-order giant component is at $p = \tfrac{1}{\binom{n}{r-s}\left(\binom{r}{s}-1 \right)}$ for each $1 \leq s \leq r-1$.
 We confirm this by using a different approach.

Throughout this paper, we will use the following asymptotic notation. For two functions $f(x)$ and $g(x)$ taking nonnegative values, we say $f(x)=O(g(x))$ (or $g(x)=\Omega(f(x))$) if $f(x) \leq C g(x)$ for some positive constant $C$ and $f(x)=o(g(x))$ if $\underset{x \rightarrow \infty}{\lim}\tfrac{f(x)}{g(x)}=0$. If $f(x)=O(g(x))$ and $g(x)=O(f(x))$, then we write $f(x)=\Theta(g(x))$.  We say an event $\X_n$  occurs with high probability (w.h.p for short)  if  the probability that $\X_n$ holds goes to $1$ as $n$ approaches infinity.

We will prove the following main theorem.

\begin{theorem}\label{main}
Let $H^r(n,p)$ be the random $r$-uniform hypergraph and $\epsilon$ be a small positive constant.
  For $1\leq s\leq r-1$,
 the following statements
  hold.
  \begin{description}
  \item[1. Subcritical phase:]
 If $p<\frac{1-\epsilon}{\left({r\choose s}-1\right){n\choose r-s }}$, then
w.h.p all $s$-th-order connected components in $H^r(n,p)$
have size $O(\ln n)$.
  \item[2. Supercritical phase:]
 If $p>\frac{1+\epsilon}{\left({r\choose s}-1\right) {n\choose r-s}}$, then
w.h.p  there is a unique giant $s$-th-order connected component
of size $\Theta({n\choose s})$ and all other $s$-th-order connected components
have size $O(\ln n)$.
Moreover, if $p=\frac{c}{ {n\choose r-s}}$ with
$c>\frac{1+\epsilon}{{r\choose s}-1}$, then  the size of the giant $s$-th-order connected component is
$(z+o(1)){n\choose s}$, where
$$z=1-\sum_{j=0}^\infty \frac{\left({r\choose
    s}j -j+1 \right)^{j-1}}{j!}c^je^{-c\left({r\choose s}j -j+1\right)}$$
     is the unique positive solution to the following equation
\begin{equation} \label{eq:gcc}
1-x=e^{c\big((1-x)^{{r\choose s}-1}-1\big)}.
\end{equation}
  \end{description}
\end{theorem}

The case $s=1$ corresponds to the  vertex-to-vertex connection.
The threshold of the emergence of the giant connected component is
$p=\frac{1}{(r-1){n\choose r-1}}$. Equivalently, the number of
edges at this critical point is about $\frac{n}{r(r-1)}$.
This is consistent to Schmidt-Pruzan and Shamir's result  \cite{PS} on $H^r(n,m)$. Similar results have been  independently proved by Cooley, Kang, and
 Koch \cite{CKK} (also see \cite{CKP}) with more careful analysis
 near the threshold. The explicit formula for the size of the giant
 component is novel.

Our method is similar to the classical approach for studying the phase transition in  Erd\H{o}s--R\'enyi
 random graph model. We will need two special Gatson-Watson branching processes-the multi-fold Poisson branching process and the multi-fold binomial branching process.
Details are
explained in Section 2. In Section 3, we will prove the main theorem.
In Section 4, we will prove the explicit formula for the size of the giant component.

\section{The Gatson-Walton branching process}
Let $\cal Z$ be a distribution over  nonnegative integers. Let $Z_1,Z_2,\ldots,$ be a countable sequence of independent identically  distributed random variables,
each having the same distribution  $\cal Z$.
A Galton-Watson process is a stochastic process $\{Y_t\}_{t=0}^{\infty}$ which
evolves according to the recursive formula
\begin{equation}
  \label{eq:gw}
Y_{t} = Y_{t-1} -1 + Z_t,  \quad \mbox{ for } t\geq 1,
\end{equation}
where $Y_0 = 1$. The process can be interpreted as follows. We think the children being born in a Depth-first Search manner. Starting from a single root
node, we can call her Eve. Eve has $Z$ children and her children are kept in a stack
(i.e., last-in first-out). Now Eve's last child
has $Z$ children and all of them are also stacked in the last-in first-out order.  A node is {\em
  dead} if its children have been exposed, and is {\em live} otherwise.  Each step, we explore  Eve's last live descendant and its children are added to the stack of live descendants.
Note that $Y_t$ is the number of live descendants  after $t$ nodes have
been explored. Equation \eqref{eq:gw} shows the recursive formula for $Y_t$.
Let $T$ be the  total number of nodes (including Eve herself) created
in the process. If the process goes on forever, then we write
$T=\infty$.

For a fixed integer $m \geq 1$, an {\em $m$-fold Poisson branching
  process}, denoted by $T^{po}_{m,c}$,
is a special Galton-Watson process with the ensemble $Z_t$
$$\Pr(Z_t=k)=
\begin{cases}
e^{-c}\frac{c^{k/m}}{(k/m)!}   & \mbox{ if } m|k  \\
0 & \mbox{ otherwise.}
\end{cases}
$$
Here $c$ is a positive constant. Each $Z_t$ is called an {\it $m$-fold Poisson random variable}. The random variable $Z_t$ can be interpreted
as the Poisson distribution $\Po(c)$ duplicated $m$ times.
Note the generating function for $Z_t$ is
$$f(x)=\sum_{\underset{m|k}{k=0}}^\infty e^{-c}\frac{c^{k/m}}{(k/m)!} x^k=
\sum_{n=0}^\infty e^{-c}\frac{c^{n}}{n!} x^{nm}=e^{c(x^m-1)}.$$
We have the following lemma.

\begin{lemma} \label{l:po}
Let $T^{po}_{m,c}$ be the $m$-fold Poisson branching
  process defined above. We have
  \begin{enumerate}
  \item If $c\leq \frac{1}{m}$, then the extinction probability is $1$.
\item If $c>\frac{1}{m}$, then the extinction probability is in $(0,1)$  and satisfies the equation
$$x=e^{c(x^m-1)}.$$
  \end{enumerate}
\end{lemma}
The proof of this lemma goes the same lines as those for proving Lemma 11.4.1 in  the monograph \cite{AS} and is omitted here.

We will  need another Galton-Watson branching process-the $m$-fold binomial branching process $T^{bin,m}_{n,p}$, where we have a countable  collection of independent identically distributed variables $Z_1,Z_2,\ldots,$ such that
 \[
\Pr(Z_i=k)=
\begin{cases}
\binom{n}{k/m} p^{k/m} (1-p)^{n-k/m} & \textrm{ if } m | k;\\
0 & \textrm{ otherwise.}
  \end{cases}
\]
Here $Z_i$ can be understood as the binomial distribution $\textrm{Bin}(n,p)$ duplicated $m$ times.

To estimate the population in the two  branching processes mentioned above,
we need the following lemma from \cite{CL3}.
\begin{lemma} \label{l:con}
  Suppose that $X_1,X_2,\ldots,X_n$ are nonnegative independent random variables.
  Let $X=\sum_{i=1}^n X_i$. Then we have the following bound on the lower tail:
  \begin{equation}
    \label{eq:lowertail}
\Pr(X \leq \E(X) -\lambda) \leq e^{-\frac{\lambda^2}{2\sum_{i=1}^n \E(X_i^2)}}.
  \end{equation}
  If further $X_i\leq M$ for all $1\leq i\leq n$,  then we have the following bound on the upper tail:
   \begin{equation}
    \label{eq:uppertail}
\Pr(X \geq \E(X) +\lambda) \leq e^{-\frac{\lambda^2}{2(\sum_{i=1}^n \E(X_i^2)+M\lambda/3)}}.
\end{equation}
\end{lemma}

We have the following lemma.
\begin{lemma}\label{l:super}
For fixed $m$ and  $c>\tfrac{1}{m}$, let $T$
be the total population of the $m$-fold Poisson branching
process $T^{po}_{m,c}$. Then there exists a constant $C=C(m,c)$ satisfying
for any positive integer $K$
$$\Pr(T=K)<e^{-CK}.$$
\end{lemma}
\begin{proof}
 Let $\{Z_t\}_{t=1}^{\infty}$ be a sequence of independent identically distributed random variables.  Each  of them have  the $m$-fold Poisson distribution. We observe
\[
\Pr(T=K)\leq \Pr(Z_1+Z_2+\ldots+Z_K=K-1).
\]
Let $X=\sum_{i=1}^KZ_i$. Then $\E(X)=\sum_{i=1}^K \E(Z_i)=mc K$
and $\sum_{i=1}^K \E(Z_i^2)=m^2(c+c^2)K$.
Applying Lemma \ref{l:con} \eqref{eq:lowertail}, we have
\begin{align*}
  \Pr(X=K-1) &\leq \Pr(X-\E(X)=-(mc-1)(K-1))\\
&< \Pr(X-\E(X)< -(mc-1)K)\\
&\leq e^{-\frac{(mc-1)^2K^2}{2m^2(c+c^2)K}}\\
&=e^{-CK}.
\end{align*}
Here $C=\frac{(mc-1)^2}{2m^2(c+c^2)}$.
\end{proof}

We recall the following Chernoff inequality which will be repeatedly used later.
\begin{theorem}{\cite{chernoff}}
\label{t:chernoff}
 Let $X_1,\ldots,X_n$ be independent random variables with
$$\Pr(X_i=1)=p_i, \qquad \Pr(X_i=0)=1-p_i.$$
We consider the sum $X=\sum_{i=1}^n X_i$
with expectation $\E(X)=\sum_{i=1}^n p_i$. Then we have
\begin{eqnarray*}
\mbox{(Lower tail)~~~~~~~~~~~~~~~~~}
\qquad \qquad  \Pr(X \leq \E(X)-\lambda)&\leq& e^{-\lambda^2/2\E(X)},\\
\mbox{(Upper tail)~~~~~~~~~~~~~~~~~}
\qquad \qquad
\Pr(X \geq \E(X)+\lambda)&\leq& e^{-\frac{\lambda^2}{2(\E(X) + \lambda/3)}}.
\end{eqnarray*}
\end{theorem}

\section{Proof of the main theorem}
\subsection{The graph branching process}
Let $H^r(n,p)$ be the random $r$-uniform hypergraph.  Throughout the rest of the paper, we let $m=\binom{r}{s}-1$. We will consider the
graph branching process of $s$-sets  described as follows.
We will maintain four families of $s$-sets $\D, \l, \N, \R$:
\begin{enumerate}
\item $\D$: the family of all {\em dead} $s$-sets (i.e., the $s$-sets
 whose neighborhoods have been explored.)
\item $\l$: the family of all {\em live} $s$-sets (i.e., the $s$-sets in the queue and
whose neighborhoods are ready to be explored.)
\item $\N$: the family of all {\em neutral} $s$-sets (i.e., the
  $s$-sets which have not been visited.)
\item $\R$: the family of all $r$-sets which has not been queried to be an
  edge or not.
\end{enumerate}
The following is the pseudo-code for the graph branching process.
\begin{tabbing}
mm\=mm\=mm\=mm\=mm\=\kill
{\bf Algorithm} Graph branching process:\+ \\
Initially
$\D\leftarrow \emptyset$, $\l\leftarrow\emptyset$, $\N\leftarrow
{V\choose s}$,
and $\R\leftarrow {V\choose r}$. \\
{\bf while} $\N\not=\emptyset$ \\
\>Select an arbitrary  $S_0\in \N$ and
move $S_0$ from $\N$ to $\l$.\\
\>{\bf while} $\l\not=\emptyset$\+\\
\>Select the last  $S\in \l$ and move $S$ from $\l$ to $\D$.\\
\> {\bf for each} $r$-set $F\in \R$ containing $S$\+\\
\>Delete $F$ from
$\R$ and query whether $F$ is an edge of $H^r(n,p)$.\\
\> If $F$ is an edge, then add all $s$-sets in ${F\choose s}\setminus (\D\cup\l)$\\
\> \hspace*{5mm} to the end of $\l$
with last-in first-out order.\\
{\bf end for each}\-\\
{\bf end while}\\
Write out the component consisting of all $s$-sets in $\D$
and reset $\D$ to $\emptyset$.\- \\
{\bf end while}\-\\
{\bf end algorithm}
\end{tabbing}

Many papers studying the phase transition in random graphs use  the Breadth-first Search
and treat $\l$ as a queue. Here we treat
$\l$ as a stack (in last-in first-out order). This is a variation of the
  Depth-first Search.
It is slightly different from
the Depth-first Search used by Krivelevich and Sudakov  in \cite{KS},
where the live set always forms a path. Here we still have a unique
path $\P$ from $S_0$ to the last visited $s$-set $S$ in the
Depth-first Search tree. The neighbors of every node in $\P$ are already
explored.  Thus the nodes in $\P$ are in $\D$. The live nodes (elements
in $\l$) are children of some node in $\P$.

Let $T^{gr}_{r,s,p}$  be the graph branching process described above. For each $s$-set, the number of $r$-sets
$F$ containing $S$ is at most ${n \choose r-s}$.  When each edge
$F$ is discovered, there are at most $m$ new  $s$-sets
added into $\l$. Thus, the graph branching process  $T^{gr}_{r,s,p}$  is dominated
by the $m$-fold binomial branching process $T^{bin, m}_{\binom{n}{r-s},p}$.  In other words,
for any integer $K$,
$$\Pr(T^{gr}_{r,s,p}>K)\leq \Pr\left(T^{bin, m}_{\binom{n}{r-s},p}>K \right).$$
This can be used to prove the statement for the subcritical phase.

\subsection{Subcritical phase}
\begin{proof}[Proof of Theorem 1 (part I):] Recall $m=\binom{r}{s}-1$ and $p=\tfrac{1-\epsilon}{ m\binom{n}{r-s}}$.
  For an $s$-set $S$, let $C(S)$ be the component containing $S$. We have
\[
\Pr(|C(S)| > K) \leq \Pr\left(T^{bin,m}_{\binom{n}{r-s},p}>K\right).
\]
Recall the definition of the $m$-fold binomial distribution. Let $\{Z_t\}_{t=1}^{\infty}$ be a countable sequence of independent identically distributed random variables, each having the $m$-fold binomial distribution. We have
\[
\Pr\left(T^{bin,m}_{\binom{n}{r-s},p}>K \right) \leq \Pr(Z_1+Z_2+\ldots+Z_K \geq K).
\]
We note for each $1 \leq i \leq K$,
\[
\E(Z_i)=1-\epsilon  \textrm{  and } \E(Z_i^2)=(1-\epsilon)m(1-p).
\]
Therefore, $\E(Z_1+Z_2+\ldots+Z_K )=(1-\epsilon)K$ and $\sum_{i=1}^K \E(Z_i^2)=(1-\epsilon)(1-p)mK$. Applying inequality \eqref{eq:uppertail} in Lemma \ref{l:con} , we get
\[
\Pr(Z_1+Z_2+\ldots+Z_K \geq K) \leq e^{- \tfrac{\epsilon^2 K}{2((1-\epsilon)(1-p)m+m\epsilon/3)}}.
\]
Now, if $K=K' \ln n$ for $K'$ large enough, we get
\[
\Pr(|C(S)| > K' \ln n) \leq \frac{1}{n^{9s}}.
\]
Therefore, the probability that there is some $S$ such that  $|C(S)| > K' \ln n $ is at most $1/n^{8s}$. Equivalently, w.h.p each component has size $O(\ln n)$ in this case.
\end{proof}

\subsection{Supercritical phase}
We  assume $p=\frac{c}{{n\choose r-s}}$ with a constant $c$
satisfying $c>\frac{1}{{r\choose s}-1}$.

At the beginning of the graph branching process, most edges are not
revealed. For each $S$, the number of edges
containing $S$ follows the binomial distribution $\textrm{Bin}((1+o(1)){n\choose
  r-s},p)$.
Each edge brings in ${r\choose
  s}-1$ new $s$-sets. However, there are two major obstacles:

\begin{enumerate}
\item  It is possible that there is an $r$-set $F$ containing $S$  has already been queried
  earlier. Thus we do not need to query it any more, which decreases the number
of $r$-sets containing $S$.

\item  When a new edge $F$ is found, not every $s$-set
of ${F\choose s}$ (other than $S$) is new. Some of them may be
already in $\l \cup \D$. This is a serious problem since it affects the number of
new $s$-sets  added in $\l$ when we have a new edge.
\end{enumerate}
To overcome the difficulty, for each $s$-set $S$, we define a new family of $r$-sets $\R_S$ as follows,
$$\R_S=\{F\in \R\colon S\subset F \textrm{ and } S'\not\subset F  \textrm{ for each } S'\in \l\cup \D
\}.$$

With this new definition, we modify the
graph branching process algorithm as follows. We fix a  small constant $\epsilon>0$ satisfying  $(1-\epsilon)c>\tfrac{1}{{r\choose s}-1}$.
If $|\R_S|\geq (1-\epsilon){n\choose r-s}$, then we query all $r$-sets
in $\R_S$; otherwise, we stop the search. Here is the pseudo-code
for the new graph branching process algorithm with halting.
This new algorithm will be denoted by GBPH, for short.

\begin{tabbing}
mm\=mm\=mm\=mm\=mm\=\kill
{\bf Algorithm} Graph Branching Process with Halting (GBPH):
\\
\>Select an arbitrary $S_0\in \N$ and
move $S_0$ from $\N$ to $\l$.\+\\
{\bf while} $\l\not=\emptyset$\+\\
Select the last element $S\in \l$\\
{\bf if } $|\R_S|<(1-\epsilon){n\choose r-s}$  \hspace{6cm} (Halting Condition)\\
% \>Exit the while loop.\\
\>Write out all $s$-sets in $\D\cup \l$ and terminate the algorithm.\\
{\bf else}\\
\>  Move $S$ from $\l$ to $\D$.\\
\> {\bf for each} $r$-set $F\in \R_S$ \+\\
\>Delete $F$ from
$\R$ and query whether $F$ is an edge of $H^r(n,p)$.\\
\> If $F$ is an edge, then add all $s$-sets in  ${F\choose s}\setminus (\D\cup \l)$\\
\>\hspace*{5mm} to the end of  $\l$ with last-in first-out order.\\
{\bf end for each}\-\\
{\bf end if}\-\\
{\bf end while}\\
Write out all $s$-sets in $\D\cup \l$.\-\\
{\bf end algorithm}
\end{tabbing}

Since we only expose the $r$-sets in $\R_S$,
each new edge $F$ brings in exactly $m$ new $s$-sets to $\l$.
Before the process is halted, for each $s$-set $S$ that we are going to explore
its neighborhood, we have
\begin{equation}
  \label{eq:Rs}
|\R_S|\geq (1-\epsilon){n\choose r-s}.
\end{equation}

%This condition guarantees that GBPH can be coupled with  a $m$-fold Poisson branching process.

 We say a component is {\em small} if its size is at most $O(\ln n)$;
it is {\bf large} if its size is $\Omega(\binom{n}{s}/(\ln n)^{2s})$.
When we run the algorithm GBPH, there are two possibilities. One the that the Halting condition is satisfied at some time. We will show in this case we get {\it  a large connect component.} (see Lemma \ref{l:super2}).  The other one  is that the condition  is never satisfied.
In this case, the algorithm will end up with  $\l=\emptyset$ and  output  a connected component. We have the following lemma.
\begin{lemma} \label{l:super1}
The probability that the Halting condition is never satisfied and  GBPH outputs a connected component is at most $(1+o(1))z$, where $z$ satisfies $0<z<1$ and
 \begin{equation}
   \label{eq:z}
z=e^{(1-\epsilon)c(z^m-1)}.
 \end{equation}
Moreover, with probability at least $1-1/n^{5s}$,  the size of the connected
component we get in this case   is at most $O(\ln n)$.
\end{lemma}
\begin{proof}
  Let $B$ be the event that the halting condition is never satisfied and $C$ is
  the component output by GBPH.  In this case,  GBPH dominates the $m$-fold binomial process $T^{bin,m}_{(1-\epsilon)\binom{n}{r-s},p}$ since $(1-\epsilon)\binom{n}{r-s}p=(1-\epsilon)c> \tfrac{1}{m}$ by the choice of $\epsilon$. Note that the size of component $C$ is always less than
  or equal to ${n\choose s}$. We have
  
  \[\Pr(B)=
    \Pr\left(B \mbox{ and }|C| \leq \binom{n}{s}\right) \leq \Pr\left(T \leq \binom{n}{s}\right) < \Pr(T \textrm{ is finite}). \]
We note  $ \Pr(T \textrm{ is finite})$ is the extinction probability of $T^{bin,m}_{(1-\epsilon)\binom{n}{r-s},p}$.
    Repeating the argument in Section 5.1 of \cite{jlr}, we know the extinction probability of the  process $T^{bin,m}_{(1-\epsilon)\binom{n}{r-s},p}$ converges to  the extinction probability of the  Poisson process $T^{po}_{m,(1-\epsilon)c}$  as $n$ approached to the infinity, which satisfies the Equation \eqref{eq:z}. Thus the first part is proved.

   In the  process GBPH, if the number of live $s$-sets $\l_t$ at some time $t$ is at least $K\ln n$, then the  probability that it outputs a connected component and the halting condition is never satisfied can be bounded from above by   $(1+o(1)z)^{K\ln n}=n^{K\ln (1+o(1)z)}$ by recalling the argument above.

    Since GBPH dominates $T^{bin,m}_{(1-\epsilon)\binom{n}{r-s},p}$, we have
    \begin{equation}
      \Pr(|\l_{t}|\geq K_1 \ln n \mbox{ and } B) \leq n^{K_1 \ln (1+o(1)z)}.
    \end{equation}
Suppose $|C| \geq K_2 \ln n$.  Set $t=K_2 \ln n$. We have
\begin{align*}
 \hspace*{-1cm} & \Pr(|C|\geq K_2 \ln n \mbox{ and } B) \\
                &\leq    \Pr(|\l_{t}|\geq K_1 \ln n \mbox{ and } B)  + \Pr(|C|\geq  K_2 \ln n \mbox{ and } |\l_{t}|\leq K_1 \ln n  \mbox{ and } B)\\
                &\leq    n^{K_1 \ln (1+o(1)z)} +\Pr(|C|\geq  K_2 \ln n \mbox{ and } |\l_{t}|\leq K_1 \ln n  \mbox{ and } B).
    % &\leq 2n^{-6s}.
    \end{align*}
Let $X_i=m\textrm{Bin}((1-\epsilon)\binom{n}{r-s},p)$ for $1 \leq i \leq t$, here $m=\binom{r}{s}-1$.  We note $\E(\sum_{i=1}^{t}) X_i=(1-\epsilon) c m t=(1+c_1)t$ for some small positive constant $c_1$ as the choice of $\epsilon$.  Recall $t=K_2 \ln n$. We  observe
\begin{align*}
\Pr(|\l_t| \leq K_1 \ln n) &\leq \Pr\left(\sum_{i=1}^t X_i \leq t+K_1 \ln n\right)\\
                                   &=\Pr\left(\sum_{i=1}^{t} X_i \leq  (1+c_1)t -(c_1t-K_1  \ln n)\right)\\
                                    & \leq e^{-(c_1K_2-K_1)^2 \ln n/ 2(1+c_1)K_2}
\end{align*}
The last inequality follows from the lower tail inequality in Theorem \ref{t:chernoff}.
We first choose $K_1$ large enough such that  $ n^{K_1 \ln (1+o(1)z)}\leq \tfrac{1}{n^{7s}}$. Then we choose $K_2$ to be sufficiently larger than $K_1$ such that   $e^{-(c_1K_2-K_1)^2 \ln n/ 2(1+c_1)K_2} \leq \tfrac{1}{n^{7s}}$.
We get $\Pr(|C| \geq K_2 \ln n ) \leq \tfrac{1}{n^{6s}}$. The union bound gives the second part.
\end{proof}

% Let ${\cal L}_t$ be the stack of live nodes ($s$-sets) at time $t$.
% If ${\cal L}_{t_1}={\cal L}_{t_2}$ for some $t_1<t_2$, then the family of
% all nodes added into $\cal D$
% during the time period $[t_1,t_2]$ forms an {\em completely-explored} subtree of
% the depth-first search tree produced by GBPH. 
% Lemma \ref{l:super1} can be strengthened into the following version. The proof is identical
% and will be omitted.
% \begin{lemma} \label{l:super2}
% With probability at least $1-1/n^{5s}$,  the size of any completely-explored subtree in GBPH is at most $O(\ln n)$.
% \end{lemma}

We keep running the algorithm GBPH until  the Halting condition is satisfied. We claim  w.h.p we get a large  component of size $\Omega(\binom{n}{s}/(\ln n)^{2s})$
(see Lemma \ref{l:super2}) when the Halting condition is satisfied.
If the claim  holds for $H^r(n,p_1)$, then it holds
for any $H^r(n,p_2)$ with $p_2>p_1$. Without loss of generality,
we assume $p=\tfrac{c}{\binom{n}{r-s}}$ with
\begin{equation}  \label{eq:p}
\frac{1}{{r\choose s}-1}<c<\frac{1-\eta}{{r-1\choose s-1}-1},
\end{equation}
for some small $\eta>0$.\\

\noindent
{\bf Fact 1:} In $H^{r}(n,p)$, with probability at least  $1-1/n^{5s}$,   for each $s$-set $S$, the number of edges containing $S$ is $O(\ln n)$.

The fact above follows from the upper tail of the Chernoff inequality (see Theorem \ref{t:chernoff}) easily  and the proof it omitted here.

We say an $s$-tree $T$ has a {\em non-empty intersection} if all $s$-stops of $T$ have a non-empty intersection.

\begin{lemma} \label{nobigtree}
  Suppose $p=\tfrac{c}{\binom{n}{r-s}}$ with $1/({r\choose s}-1)<c<(1-\eta)/({r-1\choose s-1}-1)$
  for some small $\eta>0$. With probability at least $1-1/n^{6s}$,
  every $s$-tree in $H^r(n,p)$ with non-empty intersection has size $O(\ln n)$.
\end{lemma}
\begin{proof}
For a fixed non-empty subset $A$ of size less than $s$,  let $a=|A|$  and $U=[n]\setminus A$.
Consider an $s$-tree $T$ so that all its $s$-stops contain $A$.
By deleting the common set $A$ from all $s$-stops of $T$, we get an $(s-a)$-tree $T'$ in
the random hypergraph $H^{r-a}(U,p)$.

By our choice of $p$, we have
$$p<\frac{(1-\eta)}{({r-1\choose s-1}-1) \binom{n}{r-s} }\leq \frac{1-\eta}{\left(\binom{r-a}{s-a}-1 \right) \binom{n-s}{r-s}}.$$
Therefore, we are in the subcritical case of $H^{r-a}(U,p)$. By the first part of Theorem  \ref{main}, we get with probability at least $1-1/n^{8s}$ each $(s-a)$-th connected component in $H^{r-a}(U,p)$ is  $O(\ln n) $. Note that the choices of $A$ is at most
$\sum_{a=1}^s{n\choose a}<n^{s+1}$.
The union bound completes the proof.
\end{proof}
\begin{lemma} \label{newset}
 Suppose $p=\tfrac{c}{\binom{n}{r-s}}$ with $1/({r\choose s}-1)<c<(1-\eta)/({r-1\choose s-1}-1)$
  for some small constant $\eta>0$.
When the algorithm GBPH outputs a  component $C$, let  $S_0$ be the first $s$-set entering  $C$ and $S$ be the $s$-set satisfying the Halting Condition.   If there is a proper subset $ \emptyset \not = T \subsetneq S$ such that $|\{S' \in C\colon S' \cap S \}| \geq c_1{\binom{n}{s-|T|}}/(\ln n)^k$, then with probability at least $1-1/n^{4s}$ there is a proper subset $T' \subsetneq T$ such that  $|\{S' \in C\colon S' \cap S =T'\}| \geq c_2{\binom{n}{s-|T'|}}/(\ln n)^{k+2}$, here $c_1$ and $c_2$ are constants.
\end{lemma}
\begin{proof}
Note that $C$ is a Depth-first Search tree (DFS tree for short). So it is a rooted tree.
Orient all the edges towards the root. This gives a unique root to each subtree of $C$.

To avoid ambiguity, we refer those $s$-stops as nodes (of $C$), Let ${\cal B}=\{S' \in C\colon S' \cap S=T \}$. Color  the nodes in ${\cal B}$  black and the rest of nodes (in $C$) white.
All black nodes in $C$ are naturally partitioned into a family of connected components (black subtrees of $C$). Note that each black subtree has a non-empty intersection $T$.
By Lemma \ref{nobigtree}, each black subtree has at most $K\ln n$ nodes for some positive constant $K$.
Let $\cal R$ be the set of roots of all black subtrees. Then we have
$$|{\cal R}|\geq \frac{|{\cal B}|}{K\ln n}.$$
Let ${\cal W}$ be the set of parent nodes of some node in $\cal R$. By Fact 1,  there are at most
$K'\ln n$ vertices in $\cal R$ sharing the same parent node for some positive constant $K'$. Thus we have
$$|{\cal W}|\geq \frac{|{\cal R}|}{K'\ln n}\geq \frac{|{\cal B}|}{K'K\ln^2 n}.$$

For each $S'\in {\cal W}$, we (arbitrarily) pick one (and only one) child $S''\in {\cal R}$.
Let $\cal M$ be a collection of pairs of $s$-sets  $\{(S', S'')\colon S'\in {\cal W}\}$, which forms
a matching in $C$.

For each $A\subset T$, we say a pair $(S',S'')$ is an $A$-jump if $A=T\cap S'$ and $T \subset S''$.
Since the number of choices of $A$ is at most $2^{|T|}\leq 2^s$, by the average argument,
there is an $A$ so that the number of $A$-jump pairs $(S',S'')\in {\cal M}$ is at least
$$\frac{1}{2^s}|{\cal M}|=\frac{1}{2^s}|{\cal W}|\geq c_1'\frac{\binom{n}{s-|T|}}{(\ln n)^{k+2}},$$
for some constant $c_1'$. By our construction, $T$ is not a subset of $S'$. Thus we have $|A|<|T|$.

Since the pairs in $\cal M$ form a matching of $C$, for each pair $(S',S'')\in \cal M$, there
is an edge $F:=F(S',S'')$ (in $H(n,p)$) containing $S'\cup S''$ in $H(n,p)$ is found when we explore the neighbors of $S'$.
All $F(S',S'')$'s are distinct by the algorithm DBPH.

Intuitively, to get so many $A$-jump pairs, it requires an even larger number of nodes in $C$ containing $A$.
Here is a rigorous statement.

{\bf Claim a:} For any fixed $A$ and $T$ where $A\subsetneq T$ and $|T|< s$,
let $d_A$ be the number of nodes in $C$ containing $A$ but not containing $T$.
Let $c'=\frac{c(r-s)!}{(r-s-|T|+|A|)!}$.
For any function $g(n)\geq  50s\ln n$,
if $$d_A<g(n)\frac{n^{|T|-|A|}}{2c'},$$
then with probability at most $n^{-5s}$  there are at least $g(n)$
$A$-jumps.\\

\noindent
{\it Proof of Claim a:} We list the $s$-sets containing $A$ but not $T$
as $S_1,\ldots,S_{d_A}$.  For each $1 \leq i  \leq d_A$, let $X_i$ be the random indicator variable for the $A$-jump occurring at $S_i$.
Let $X=\sum_{i=1}^{d_A}X_i$ be the number of $A$-jumps.
For $1 \leq i \not = j \leq d_A$, the  collection of $r$-sets we query for $S_i$ is disjoint with the one for $S_j$ by the algorithm GBPH, so $X_i$'s are independent random variables.
We note  an  $A$-jump happens at $S_i$ if and only if there is an $r$-set $F$ satisfying $T \subset F$ is found as an edge when we explore the neighbors of $S'$. Here we notice $A\subsetneq T$ and $
T \subsetneq S''$.
Therefore,  the $A$-jump happens at $S_i$ with probability at most
\begin{equation} \label{eq6}
{n-s \choose r-s-|T|+|A|}p\leq \frac{c(r-s)!}{(r-s-|T|+|A|)!n^{|T|-|A|}}
=\frac{c'}{n^{|T|-|A|}}.
\end{equation}
We have
$$\E(X)\leq d_A  \frac{c'}{n^{|T|-|A|}}< \frac{g(n)}{2}.$$
 Applying the upper tail of Chernoff's inequality (see Theorem \ref{t:chernoff}), we have
\begin{align*}
  \Pr(X\geq g(n)) &\leq \Pr\left(X-\E(X)\geq \frac{g(n)}{2}\right)\\
&\leq e^{-3g(n)/16}\\
&<\frac{1}{n^{7s}}.
\end{align*}

Choose $g(n)=\frac{c_1'}{2}\frac{\binom{n}{s-|T|}}{(\ln n)^{k+2}}$.
Since we have at least $2g(n)$ $A$-jump pairs, by Claim a,  w.h.p we have
$$d_A\geq g(n)\frac{n^{|T|-|A|}}{2c'}> \frac{c_1'}{2}\frac{\binom{n}{s-|T|}}{(\ln n)^{k+2}}
\frac{n^{|T|-|A|}}{2c'}>c_2{\binom{n}{s-|A|}}/(\ln n)^{k+2},$$
for some constant $c_2$. The proof of Lemma is finished by taking $T'=A$.
\end{proof}

\begin{lemma}\label{l:super2}
If  the algorithm GBPH outputs a  component because of the Halting condition, then  w.h.p it  contains an $s$-tree of size $\Omega\left({n\choose s}/{(\ln  n)^{2s}}\right)$.
\end{lemma}
\begin{proof}
We keep running the algorithm GBPH $t$ times.  For each  $1 \leq i \leq t$,  let $X_i$ be the event that the $i$-th  running of the algorithm GBPH gives a connected component  $C_i$ and the Halting condition is not satisfied.
Let $X$ be the event that the algorithm GBPH outputs $k$ connected components and the Halting condition is not satisfied.
We note
\[
\Pr(X)= \Pr(X_1X_2\cdots X_t) = \Pr(X_1) \Pr(X_2|X_1) \cdots \Pr(X_t|X_1 \cdots X_{t-1}).
\]
We define $V_0=\{1,\ldots,n\}$, $U_i=\cup_{S \in C_i} S$ for $1 \leq i \leq t$, and $V_i=V_0\setminus (\cup_{j=1}^{i-1} U_j)$ for each $2 \leq i \leq k$.   We note  $\Pr(X_i|X_1 \cdots X_{i-1})$  is bounded from above by the probability that  a running of the algorithm GBPH in $H^{r}(V_i, p)$ gives a connected component and the Halting condition is not satisfied..   By  Lemma \ref{l:super1},  with probability $1-1/n^{5s}$ we have $|C_i| \leq K \ln n $ for some constant $K$.  Thus, with  probability at least $1-t/n^{5s}$ we have $|V_i| \geq n - is K \ln n$ for each $1 \leq i \leq t$.  We choose $t=K' \ln n $ for some constant $K'$ and let $U=\{1,\ldots,n-sKK'\ln n\}$. Then $|V_i| \geq |U|$ for each $1 \leq i \leq t$. Thus,   $\Pr(X_i|X_1 \cdots X_{i-1})$  is bounded from above by the the probability that a running of the algorithm GBPH in $H^{r}(U, p)$ gives a connected component.  By  Lemma \ref{l:super1},  we get  $\Pr(X_i|X_1 \cdots X_{i-1}) \leq (1+o(1))z$. Now, if we take $K'$ large enough, then we get
\[
  \Pr(X)=\Pr(X_1) \Pr(X_2|X_1) \cdots \Pr(X_t|X_1 \cdots X_{t-1}) \leq ((1+o(1))z)^{K' \ln n}=
  o \left(\frac{1}{n^{6s}}\right).
\]

Let $C$ be the  component given by the algorithm GBPH and $S$ be the $s$-set such that
$|\R_S|\leq (1-\epsilon){n\choose r}.$
In other words,
\begin{equation}
  \label{eq:S}
\sum_{S' \in C}{n\choose r-|S\cup S'|}\geq \epsilon{n\choose r-s},
\end{equation}
We make convention $\binom{n}{r-|S \cup S'|}=0$ in the case of
 $|S\cup S'|> r$.  We use $\cal Q$ to denote the collection of all subsets of $S$.
For each $T \in \cal Q$, we define a constant
\[
c_T=\frac{\epsilon (r-2s+|T|)!(s-|T|)!}{(r-s)!2^{|S+1|}}.
\]
We claim there is a proper subset $T$ of  $S$ such that
\[
|\{S' \in C : S\cap S'=T\}| \geq c_T \binom{n}{s-|T|}.
\]
Otherwise, we have
\begin{align*}
\sum_{S' \in C}{n\choose r-|S\cup S'|}&=\sum_{T \in \cal Q} \sum_{\underset{S' \in C}{S' \cap S =T}} \binom{n}{r-2s+|T|} \\
                                                            &< \sum_{T \in \cal Q}  c_T  \binom{n}{r-2s+|T|}  \binom{n}{s-|T|}\\
                                                            & <  \sum_{T \in \cal Q}  \frac{\epsilon \binom{n}{r-s}}{2^{|S|}}\\
                                                            &=\epsilon \binom{n}{r-s}
\end{align*}
Contradiction to Inequality \eqref{eq:S}. The claim is proved.

If $T=\emptyset$, then we have $|C| \geq  |\{S' \in C : S\cap S'=\emptyset\}| \geq c_T \binom{n}{s}$ and we are done. Otherwise, we will apply Lemma \ref{newset} recursively.   Initially, we feed   Lemma \ref{newset} with the subset $T=T_0$ claimed above and let the output  be $T_1$.
For each  $i \geq 1$, we apply  Lemma \ref{newset} with $T_{i-1}$ and   let $T_i$ be the output.  We note  $|T_i| \leq s-i-1$  for each   $i \geq 0$.
Moreover, there is constant $c_i$ such that
\begin{equation} \label{eq:empty}
|\{S' \in C: S' \cap S =T_i\}| \geq \frac{c_i \binom{n}{s-|T_i|}}{(\ln n)^{2i} }   \textrm{  for each } i \geq 0 .
\end{equation}
Let $j$ be the smallest number such that $T_j=\emptyset$. We have $j \leq s-1$ and the lemma follows from \eqref{eq:empty}.  We note  the error probability is at most $K'\ln n/n^{5s}+s/n^{4s} \leq 1/n^{3s}$.
The lemma follows from the union bound.
\end{proof}

%After the  algorithm first stops, we query all $r$-sets which we have not checked. Finally, we output all components.
We already proved that there exists at least one large  component. We next show the large  components is unique.
\begin{lemma}\label{l:unique}
W.h.p the large  component is unique.
\end{lemma}
\begin{proof}
Let $C$ be  a  component such that $|C|=\eta \binom{n}{s}/(\ln  n)^{2s}$ for some $\eta >0$.  For $s \geq 2$, we will  prove that $C$ must have some nice properties.

 For each $U \subset [n]$ such that $1 \leq |U| \leq s-1$ , let
\[\Gamma_C(U)=\left\{S\in C \textrm{ such that } U \subset S\ \right \} \textrm{ and } d_C(U)=|\Gamma_C(U)|.\]
For each $1 \leq i \leq s-1$,  let $f_i=(s+2)^2/2-i(i+2)/2$.
 We next show the following claim.

 {\bf Claim b:} For each $1 \leq i \leq s-1$,  with  probability at least $1-1/n^{f(i)}$ we have
$\binom{[n]}{i} \subset \cup_{S \in C} \binom{S}{i}$  and there exists some  constant $\delta_i >0$ such that $ d_C(U) \geq \delta_i n^{s-i}/ (\ln  n)^{2s}$ for each $U \in \binom{[n]}{i}$.

\noindent
{\it Proof of Claim b:}
We prove the claim by induction on $i$.  For the case where $i=1$, we suppose there is some $1 \leq i \leq n$ such that $i \not \in \cup_{S \in C} S$.  For each $S \in C$, let $T$ be an arbitrary subset of $S$ with size $s-1$ and $S'=T \cup \{i\}$. Then the number of $r$-sets containing $S \cup S'$ is $\binom{n}{r-s-1}$. Let $S$ run over all $s$-sets of $C$. Note that each $r$-set counts at most $\binom{r}{s+1}$ times. The expected number of $r$-sets containing $S \cup S'$ for some $S \in C$ is at least $p\binom{n}{r-s-1}|C|/\binom{r}{s+1} = \Omega(n^{s-1}/(\ln  n)^{2s})$. As we assume $s \geq 2$, the lower tail of Chernoff's inequality (see Theorem \ref{t:chernoff}) shows that with probability at least $1-1/n^{f(1)+1}$ there are at least $\Omega(n^{s-1}/(\ln  n)^{2s})$ $r$-sets containing $S \cup S'$ for some $S \in C$, i.e.,
$i \in \cup_{S \in C} S$ and $d_C(i) \geq \delta_1 n^{s-1}/(\ln  n)^{2s}$ for some positive constant $\delta_1$.
 Since we have at most $n$ choices for $i$, the base case of the claim follows from the union bound.

 For the inductive step,  suppose there is some $U \in \binom{[n]}{i}$ such that
$U \not \in \cup_{S \in C } \binom{S}{i}$. Let $U'$ be an arbitrary subset of $U$ with size $i-1$. By the inductive hypothesis, with probability at least $1-1/n^{f(i-1)}$ we have $d_C(U') \geq \delta_{i-1}n^{s-i+1}/ (\ln  n)^{2s}$.  For each $S \in \Gamma_C(U')$, let $T$ be a subset of $S$ such that $|T|=s-1$ and  $U' \subset T$.   We define a new $s$-set $S'=U \cup T$. By the same argument as in the base case, the expected number of $r$-sets containing $S \cup S'$ for some $S \in C$ is at least $\Omega(n^{s-i}/(\ln n)^{2s})$. The lower tail of Chernoff's inequality (see Theorem \ref{t:chernoff}) shows with probability at least $1-1/n^{f(i-1)+0.5}$ there are at least $\Omega(n^{s-i}/(\ln  n)^{2s})$ $r$-sets containing $S \cup S'$ for some $S \in C$,
 $U \in \cup_{S \in C} \binom{S}{i}$ and $d_C(U) \geq \delta_i n^{s-i}/(\ln  n)^{2s}$ for some positive constant $\delta_i$. As we have at most $n^i$ choices for $U$, the inductive step follows from the union bound.

We are now ready to prove the lemma. We first consider the case of $s \geq 2$.
Suppose there are two $s$-sets $S_1$ and $S_2$ which are contained in components  $C_1$ and $C_2$  respectively, here $C_1$ and $C_2$ have size at least $\eta \binom{n}{s}/\ln^{2s} n$ for some $\eta >0$.  Let $U$ be an arbitrary $(s-1)$-subset of $[n]$. By Claim b, with probability at least $1-2/n^{2s+2}$ we have $d_{C_i}(U) \geq \delta_{s-1}n /(\ln  n)^{2s}$ for each $i \in \{1,2\}$.  For  each $S_1' \in \Gamma_{C_1}(U)$ and $S_2' \in \Gamma_{C_2}(U)$, the number of $r$-sets containing $S_1' \cup S_2'$ is $\binom{n}{r-s-1}$. With the lower bounds on $d_{C_1}(U)$ and $ d_{C_2}(U)$, the expected number of $r$-sets containing $S_1' \cup S_2'$ for some $S_1' \in \Gamma_{C_1}(U)$ and $S_2' \in \Gamma_{C_2}(U)$ is $\Omega\left(n^2/(\ln  n)^{4s}\right)$, here we note each $r$-set counts at most $r^2$ times.  Now the lower tail of Chernoff's inequality (see Theorem \ref{t:chernoff}) gives with  probability at least $1-1/n^{2s+2}$ $S_1'$ is adjacent to $S_2'$ for some $S_1' \in \Gamma_{C_1}(U)$ and $S_2' \in \Gamma_{C_2}(U)$. Therefore, with  probability at least $1-3/n^{2s+2}$ we have $S_1$ and $S_2$ are in the same component. As there are at most $n^{2s}$ choices for the two $s$-set $S_1$ and $S_2$, we get w.h.p the large component is unique. We can prove the case where $s=1$ exactly by the same lines as above.
\end{proof}

\begin{proof}[Proof of Theorem \ref{main} (part II):]
By Lemma \ref{l:super2} and Lemma \ref{l:unique}, the giant
$s$-th-order component exists and is unique. We need to show that the
size of the  giant $s$-th-order component concentrates on
$(z+o(1))\binom{n}{s}$, where
$z$ satisfies equation \eqref{eq:gcc}.

Recall $m=\binom{r}{s}-1$ and $p=\frac{c}{\binom{n}{r-s}}$, where $c>\frac{1+\epsilon}{\binom{r}{s}-1}$.
  Let $T^{{bin},m}_{n,p}$ denote the $m$-fold binomial branching process.
For an $s$-set $S$, let $C(S)$ be the $s$-th-order connected component
containing $S$.
For a fixed $t$, we observe
\begin{equation} \label{eq:size}
\Pr\left(T_{ \binom{n-s-st}{r-s}, p}^{bin,m} \geq t \right) \leq  \Pr(
|C(S)| \textrm{ is not small} ) \leq \Pr\left(T_{ \binom{n}{r-s}, p}^{bin,m} \geq t \right).
\end{equation}
In  the graph branching process, when we explore the neighbors of the current $s$-set $S$, we query all $r$-sets (have not checked) containing $S$.  To get the lower bound  in \eqref{eq:size}, we can consider only the $r$-set $F$ such that
\[S \subset F \textrm{ and }  F \cap \left( \cup_{S' \in \l \cup \D} (S' \setminus S) \right)=\emptyset.
\]
Note that this kind of $F$ brings in $\binom{r}{s}-1$ new $s$-sets if
$F$ is an edge in $H^r(n,p)$. We observe the number of $F$ satisfying
the condition above is at least $\binom{n-s-st}{r-s}$. Then we get the
lower bound. The upper bound follows similarly.

Intuitively, $\Pr\left(T_{ \binom{n-s-st}{r-s}, p}^{bin,m} \geq t
\right)$ and  $\Pr\left(T_{ \binom{n}{r-s}, p}^{bin,m} \geq t \right)$
are approximately $\Pr(T_{m,c}^{po} \geq t)$, which is arbitrarily
close to $\Pr(T_{m,c}^{po}=\infty)$ as both $t$ and $n$ approach infinity.

Here is a rigorous argument. Let $z=\Pr(T_{m,c}^{po}=\infty)$.
By Lemma \ref{l:po}, $z$ is the unique positive root of equation \eqref{eq:gcc}.

For any $\varepsilon>0$,  we can find an
integer $t=t(m,c,\varepsilon)$ so that $\sum_{K\geq t}^\infty
e^{-CK}<\varepsilon$, where $C$ is the constant given by Lemma \ref{l:super}.
Thus, we have
$$z<\Pr(T^{po}_{m,c}\geq t)\leq
z +\varepsilon.$$
We notice that $t$ is independent of $n$. Since the Poisson
distribution  is the limit of Binomial distribution, we have
\begin{align*}
\lim_{n\to \infty}\Pr\left(T_{ \binom{n}{r-s}, p}^{bin,m} \geq t
\right)
&=1-\sum_{i=0}^{t-1} \lim_{n\to \infty}\Pr\left(T_{ \binom{n}{r-s},
    p}^{bin,m} =i\right) \\
&=1- \sum_{i=0}^{t-1} \Pr\left(T^{po}_{m,c} =i\right) \\
&= \Pr\left(T^{po}_{m,c} \geq t \right) \\
&< z+\varepsilon.
\end{align*}
Similarly, we have
$$\lim_{n\to \infty}\Pr\left(T_{ \binom{n-s-st}{r-s}, p}^{bin,m} \geq t \right)
=\Pr\left(T^{po}_{m,c} \geq t\right) >z.$$
There exists a constant $n_0=n_0(m,c,\varepsilon)$ such that for
$n\geq n_0$ we have
\begin{align}
\label{eq:sizeub}
\Pr\left(T_{ \binom{n}{r-s}, p}^{bin,m} \geq t
\right)&<z+2\varepsilon\\
\label{eq:sizelb}
\Pr\left(T_{ \binom{n-s-st}{r-s}, p}^{bin,m} \geq t
\right)&>z-\varepsilon.
\end{align}
Combining Equations \eqref{eq:size}, \eqref{eq:sizeub}, and
\eqref{eq:sizelb}, we have
$$z-\varepsilon<  \Pr(
|C(S)| \textrm{ is not small} ) <z+2\varepsilon.$$
Since $\varepsilon$ is arbitrary, we conclude that $\Pr(
|C(S)| \textrm{ is not small} )$ concentrates around $z$.

 Now, each $s$-set has probability $z+o(1)$ to be in the large connected
component. The expected number of $s$-sets in the large connected
component is $(z+o(1))\binom{n}{s}$.
Since we already showed the large
connected component is unique, this shows that the size of the giant
$s$-th-order connected component is $(z+o(1))\binom{n}{s}$.

We will postpone the proof of the explicit formula for  $z$ until next section.
\end{proof}

\section{The explicit formula for the size of the giant component}
In this section, we will study the number of small components of
$H^r(n,p)$. We will use it to deduce an explicit formula for  the size
of the giant component
 in the supercritical phase.

We have the following theorem.
\begin{theorem}\label{t2}
 For any $r\geq 2$, $1\leq s <r$, any constant $c>0$, and
 sufficiently large $n$,  let $p=\frac{c}{{n\choose r-s}}$ be the
 probability. For any positive integer $j$ and any $s$-set $S$, the
 probability that $S$ is in an $s$-th-order component of exactly $j$ edges in $H^r(n,p)$
 is
$$(1+o(1))\frac{({r\choose s}j -j+1)^{j-1}}{j!}c^je^{-c({r\choose s}j -j+1)}.$$
\end{theorem}

Note an $s$-set $S$ is in the giant component if and only if $S$ is
not in any small component. As a corollary, we have
$$z=1- \sum_{j=0}^\infty\frac{({r\choose s}j -j+1)^{j-1}}{j!}c^je^{-c({r\choose s}j -j+1)}.$$
Thus we proved  the explicit formula for $z$ in Theorem \ref{main}.

\begin{proof}
Typical $s$-th-order components of exactly $j$ edges form a special class of
$s$-trees,
called {\em $(r,s)$-trees}.
An $r$-uniform hypergraph is called an
 {\em $(r, s)$-tree} if it can be obtained from a single edge by consecutively adding edges so that every new edge contains $r-s$ new vertices while
its remaining $s$ vertices are covered by an already existing edge.
Pikhurko \cite{Pikhurko} proved that the number $b_j$ of distinct vertex
labelled $(r,s)$-trees with $j$ edges is given by
$$b_j=\frac{(j(r-s)+s)!(j{r\choose s}-j+1)^{j-2}}{j!s!((r-s)!)^j}.$$

For a fixed $(r,s)$-tree $T$ with $j$ edges (as a subgraph of
$K^r_n$), let $f(T)$ be the number of edges in $K^r_n$ which is not an
edge of $T$ but intersects at least one edge of $T$ on $s$ vertices or
more. Since $T$ has exactly $j{r\choose s}-j+1$ $s$-sets, we have
\begin{equation}
  \label{eq:fT}
\left(j{r\choose s}-j+1\right){n-j(r-s)-s\choose r-s}\leq
f(T)\leq \left(j{r\choose s}-j+1\right){n\choose r-s}.
\end{equation}
The probability that $T$
forms a $s$-th order component of $H^r(n,p)$ is given by
$$p^j(1-p)^{f(T)}=(1+o(1))p^je^{-(j\binom{r}{s}-j+1)p\binom{n}{r-s}}  = (1+o(1))p^j e^{-\left(j{r\choose s}-j+1\right)c}.
$$
Therefore, the expected number of components in $H^{r}(n,p)$ forming an  $(r,s)$-trees
with $j$ edges is
\begin{align*}
& (1+o(1)) {n\choose s+(r-s)j}b_jp^j e^{-\left(j{r\choose
       s}-j+1\right)c} \\
&=(1+o(1))
\frac{n^{s+(r-s)j}}{(s+(r-s)j)!}\frac{(j(r-s)+s)!(j{r\choose
    s}-j+1)^{j-2}}{j!s!((r-s)!)^j}p^j e^{-\left(j{r\choose
       s}-j+1\right)c}\\
&=(1+o(1))\frac{n^s}{s!} \frac{(j{r\choose
    s}-j+1)^{j-2}}{j!}
\left(\frac{n^{r-s}}{(r-s)!}p\right)^je^{-\left(j{r\choose
       s}-j+1\right)c} \\
&=(1+o(1)){n\choose s} \frac{(j{r\choose
    s}-j+1)^{j-2}}{j!} c^je^{-\left(j{r\choose
       s}-j+1\right)c}.
\end{align*}
 Note that a component with $j$ edges not forming an $(r,s)$-tree has
less vertices than an $(r,s)$-tree with $j$ edges.  Since we view $j$ as a constant,  the expected number of such
components is a lower-order term.
Each $(r,s)$-tree with $j$ edges has exactly $j{r\choose
    s}-j+1$ $s$-sets and the number of $s$-sets in total is $\binom{n}{s}$.
Therefore,  for a fixed $s$-set $S$, the probability $S$ is
  in an component of $j$ edges is
$$(1+o(1)) \frac{(j{r\choose
    s}-j+1)^{j-1}}{j!} c^je^{-\left(j{r\choose
       s}-j+1\right)c}.$$
 The proof of Theorem \ref{t2} is finished.
\end{proof}

\noindent
{\bf Acknowledgment:}  Authors would like to thank two anonymous referees for their valuable comments which greatly improved the presentation of the paper.

\end{document}